\documentclass[11pt]{amsart}

\usepackage{eucal,amsmath,amsthm,amsfonts,verbatim,amsopn,amssymb, amscd}

\usepackage[all]{xy}

\newcommand\dmo{\DeclareMathOperator}

\newtheorem{theorem}{Theorem}[section]
\newtheorem*{nthm}{Theorem}
\newtheorem{conjecture}[theorem]{Conjecture}

\newtheorem{lemma}[theorem]{Lemma}
\newtheorem{corollary}[theorem]{Corollary}
\newtheorem{proposition}[theorem]{Proposition}
\newtheorem{comm}{(COMMENT)}

\theoremstyle{definition}

\newtheorem{definition}[theorem]{Definition}
\newtheorem{notation}[theorem]{Notation}
\newtheorem{example}[theorem]{Example}
\newtheorem{remark}[theorem]{Remark}
\newtheorem{remarks}[theorem]{Remarks}

\newtheorem{question}[theorem]{Question}
\newcommand\bthm{\begin{theorem}}
\newcommand\ethm{\end{theorem}}
\newcommand\bcn{\begin{conjecture}}
\newcommand\ecn{\end{conjecture}}
\newcommand\bla{\begin{lemma}}
\newcommand\ela{\end{lemma}}
\newcommand\bco{\begin{corollary}}
\newcommand\eco{\end{corollary}}
\newcommand\bpro{\begin{proposition}}
\newcommand\epro{\end{proposition}}
\newcommand\bcm{\begin{comm}}
\newcommand\ecm{\end{comm}}
\newcommand\bdf{\begin{definition}}
\newcommand\edf{\end{definition}}
\newcommand\bnot{\begin{notation}}
\newcommand\enot{\end{notation}}
\newcommand\bex{\begin{example}}
\newcommand\eex{\end{example}}
\newcommand\brm{\begin{remark}}
\newcommand\erm{\end{remark}}
\newcommand\brms{\begin{remarks}}
\newcommand\erms{\end{remarks}}
\newcommand\bqu{\begin{question}}
\newcommand\equ{\end{question}}

\newcommand\bea{\begin{array}{r@{,\ \ \ \ \ \ \ }l}}
\newcommand\ena{\end{array}}
\newcommand\barr{\begin{array}}
\newcommand\earr{\end{array}}
\newcommand\beq{\begin{equation}}
\newcommand\eeq{\end{equation}}
\newcommand\bal{\begin{align}}
\newcommand\eal{\end{align}}

\newcommand\bqa{\begin{eqnarray*}}
\newcommand\eqa{\end{eqnarray*}}
\newcommand\bqan{\begin{eqnarray}}
\newcommand\eqan{\end{eqnarray}}
\newcommand\bit{\begin{itemize}}
\newcommand\ben{\begin{enumerate}}
\newcommand\een{\end{enumerate}}
\newcommand\eit{\end{itemize}}
\newcommand\bfi{\begin{figure}[htbp]}
\newcommand\efi{\end{figure}}
\newcommand\bce{\begin{center}}
\newcommand\ece{\end{center}}
\newcommand\bpr{\begin{proof}}
\newcommand\epr{\end{proof}}

\newcommand{\ignore}[1]{}

\newcommand\noi{\noindent}

\newcommand\msk{\medskip}

\newcommand\sub{\subseteq}

\newcommand\ra{\rightarrow}

\newcommand\abar{\overline{a}}
\newcommand\bbar{\overline{b}}
\newcommand\cbar{\overline{c}}
\newcommand\ebar{\overline{e}}

\newcommand\sbar{\overline{s}}
\newcommand\xbar{\overline{x}}
\newcommand\ybar{\overline{y}}
\newcommand\zbar{\overline{z}}

\newcommand\mc{\mathcal}

\newcommand\mcO{\mathcal{O}}

\newcommand\FF{\mathbb{F}}
\newcommand\QQ{\mathbb{Q}}
\newcommand\ZZ{\mathbb{Z}}
\newcommand\NN{\mathbb{N}}

\newcommand\mfa{\mathfrak{a}}
\newcommand\mfp{\mathfrak{p}}
\newcommand\mfq{\mathfrak{q}}
\newcommand\mfm{\mathfrak{m}}

\newcommand\pd{\partial}

\newcommand\dl{\delta}

\dmo{\Spec}{Spec}
\dmo{\Specd}{Spec_D}
\dmo{\Spf}{Spf}
\dmo{\Art}{Art}
\dmo{\Ext}{Ext}
\dmo{\Hom}{Hom}
\dmo{\Def}{Def}
\dmo{\Av}{AV}
\dmo{\Avbt}{AVBT}

\begin{document}

\title{A differential Chevalley theorem}
\author{Eric Rosen}
\address{Department of Mathematics\\
Massachusetts Institute of Technology\\
Cambridge, MA 02139}
\email{rosen@math.mit.edu}
\urladdr{http://math.mit.edu/~rosen}
\thanks{I would like to thank Victor Kac for asking me the question that led to the main
result in this paper.  I am also grateful to Ehud Hrushovski for providing a proof of
Proposition~\ref{pro:Hr}.}

\date{\today}

\begin{abstract}
We prove a differential analog of a theorem of Chevalley on
extending homomorphisms for rings with commuting derivations,
generalizing a theorem of Kac.  As a corollary, we establish that,
under suitable hypotheses, the image of a differential scheme
under a finite morphism is a constructible set.  We also obtain
a new algebraic characterization of differentially closed fields.
We show that similar results hold for differentially closed fields
that are saturated, in the sense of model theory.  
In characteristic $p > 0$, we obtain related results and 
establish a differential Nullstellensatz.  
Analogs of these theorems for difference fields are also considered.
\end{abstract}

\maketitle

\section{Introduction}

Chevalley proved the following theorem about extensions of
homomorphisms.
\begin{nthm}
Let $S$ be a noetherian integral domain and $R$ a subring such
that $S$ is finitely generated over $R$.  For any nonzero $b \in S$,
there is a nonzero $a \in R$ such that any homomorphism 
to an algebraically closed field, $\phi: R \ra K$, with $\phi(a) \neq 0$, 
lifts to a homomorphism $\psi: S \ra K$, with $\psi(b) \neq 0$.
\end{nthm}
\noindent
This implies a basic fact of algebraic geometry.
Let $f: X \ra Y$ be a morphism of finite type of noetherian schemes.
Then $f(X)$ is a constructible set.  This is closely related to
Tarski's elimination of quantifiers theorem for algebraically closed 
fields.

Recently Kac~\cite{Kac} established an analog of Chevalley's theorem in 
differential algebra, for rings equipped with a single derivation.
Our main result generalizes Kac's theorem to rings with finitely 
many commuting derivations.  As a corollary, we obtain a geometric Chevalley 
theorem for differential schemes with commuting derivations, extending a 
result of Buium.   We also provide the following new characterization of 
differentially closed fields, suggested by another result of Kac.
A differential field $K$ is differentially closed if and only if for
any finitely generated differential $K$-algebra $S$ with no zero divisors
and any nonzero $b \in S$, there is a homomorphism $\psi: S \ra K$ 
with $\psi(b) \neq 0$.

The proof of our main result follows the general strategy used by Kac,
but makes essential use of the quantifier elimination theorem for
differentially closed fields with commuting derivations due to
McGrail~\cite{TM} (see also~\cite{Yaf}) which is closely related to the elimination theorem
of Seidenberg~\cite{AS}.  As with Tarski's theorem above, 
our differential Chevalley theorem is related to quantifier elimination,
though it is certainly not equivalent.  In particular, while we show, following 
Kac, that the conditions of our main result provide a characterization of differentially
closed fields, there are many other differential fields with quantifier 
elimination~\cite{HI}.
(This result of Hrushovski and Itai stands in contrast to the fact
that the infinite fields with elimination of quantifiers are exactly those
that are algebraically closed.)

In Section~\ref{sec:sat}, we consider variations on the earlier
results where the differential algebra $S$ is not finitely generated
over the subring $R$.  Assuming that $S$ is integral over $R$, we obtain
analogs of previous theorems with the additional condition that the 
differentially closed field under consideration is sufficiently 
{\em saturated} in the sense of model theory.  (Saturation is a certain
kind of largeness property which can be defined, in this context, in
terms of the simultaneous solution of infinite sets of differential
polynomial equations and inequations over a differential field.)
The results make use of a differential version of the going up theorem
from commutative algebra, which may be of independent interest.

Differential fields in characteristic $p$ are considered in Section~\ref{sec:char.p}.
We give an easy example showing that the analog of the main theorem in this
context fails.  We then establish a related statement, which also yields a 
new characterization of the differentially closed fields in positive characteristic.
As a corollary, we obtain a differential Nullstellensatz.

In Section~\ref{sec:ACFA}, we consider difference fields, that is, fields equipped
with a distinguished isomorphism.  Although there are many similarities between
differential and difference algebra, we show that most of our results become
false in the context of difference fields.  In a separate paper~\cite{Ros}, though,
we establish some related results for difference fields, including a 
Nullstellensatz, using some ideas developed here.

In the final section, we give a new proof of Chevalley's original theorem
along the lines of the argument of our main theorem, using Tarski's quantifier
elimination for algebraically closed fields.  We also consider whether there is an
abstract model theoretic version of Chevalley's theorem, and pose
an open question in this direction.

\section{Background}

\subsection{Preliminaries}

Throughout this paper, we will consider rings equipped with $N$ commuting
derivations $\pd_1, \ldots , \pd_N$, for some fixed number $N$, which we 
also assume to contain $\QQ$.
We let $\Delta = \{\pd_1, \ldots , \pd_N\}$, and call such rings $\Delta$-rings,
or simply differential rings.  In the literature they are also referred to as Ritt algebras.
Let $\Theta$ be the set of formal expressions  $\pd_1^{e_1} \cdots \pd_N^{e_N}$,
$e_i \in \NN$.  Given a differential ring $R$, the {\em differential polynomial ring}
in $m$ variables, denoted $R\{x_1, \ldots x_m\}$, is the polynomial ring over $R$
in the variables $\theta x_i$, $\theta \in \Theta, i \leq m$, made into a  
$\Delta$-ring in the obvious way.  We call each expression $\theta x_i$ a 
{\em differential variable}.  

Given a $\Delta$-ring $R$, a {\em differential ideal} $I \sub R$ is an ideal
such that for all $r \in I$ and $i \leq N$, $\pd_i(r) \in I$.  In this case,
the quotient $R/I$ is also a $\Delta$-ring.  For any set $A \sub R$, let
$[A]$ denote the differential ideal generated by $A$.  
Given an ideal $I \sub R$, and an element $a \in R$, then $I : a^\infty$  
is by definition the ideal $\{b \in R : a^nb \in I \mbox{ for some } n\}$.

A homomorphism between $\Delta$-rings is a homomorphism that commutes 
with each of the derivations.
A {\em differential algebra} $S$ over a differential field $K$ is a differential ring
together with an embedding of $K$ into $S$.
For $R$ a differential ring, $S \sub R$ a differential subring, and $t \in R$,
then $S\{t\}$ denotes the differential subring generated by $S$ and $t$.
We say that $t$ is {\em differentially transcendental} over $S$ if for every nonzero
polynomial $f(x) \in S\{x\}$, $S(t) \ne 0$.  Otherwise, $t$ is differentially algebraic.

\subsection*{Conventions and notation}

We will often write polynomial for differential polynomial, and otherwise 
drop the word ``differential'' when our meaning is clear from context.
Generally, we will reserve the letters $x,y,z$ for variables in polynomial rings.
Let $R$ be a differential ring, $R\{y\}$ the differential polynomial ring
in one variable.  Given a polynomial $f(y)$ in $R\{y\}$, we 
will sometimes write $f(y)$ as $\hat{f}(c_0, ..., c_n, y)$ to indicate
that the coefficients that occur in $f(y)$ are among the set $\{c_0, \ldots , c_n\}$.
For example, given $g(y) = 5x^2 + rx + s$, with $r,s \in R$,
then $\hat{g}(r,s,y) = 5x^2 + rx +s$.

Given a polynomial $\hat{f}(c_0, \ldots c_n,y)$, we can replace each
occurence of a coefficient $c_i$ by a variable $z_i$, to obtain a polynomial,
$\hat{f}(z_0, \ldots , z_n, y) \in \ZZ\{z_0, \ldots, z_n, y\}$, a differential
polynomial ring in $n+2$ variables.  Thus, in the above example, 
$\hat{g}(z_0,z_1,y) = 5x^2+z_0x+z_1$, which is in $\ZZ\{z_0,z_1,x\}$.
Also, given a polynomial $\hat{f}(c_0, \ldots, c_n,y)$ and elements 
$e_0, \ldots , e_n$ in $R$, we may write $\hat{f}(e_0, \ldots , e_n, y)$ for
the polynomial in $R\{y\}$ that one obtains by replacing each occurence
of $c_i$ by $e_i$.

For brevity, we will usually write, e.g., $\cbar$ for $c_0, \ldots , c_n$, likewise, $\zbar$, etc.

\subsection{Some differential algebra}

We summarize some information about prime differential ideals in differential polynomial
rings that we will need later.  An exhaustive reference for this material is Kolchin's
book~\cite{EK1}.

\bdf
\label{df:sep}
Let $R$ be a $\Delta$-ring, $R\{y_1, \ldots y_m\}$ the $\Delta$-polynomial ring in 
$m$ variables.  To each differential variable $\pd_1^{e_1}\cdots \pd_N^{e_N}y_j$,
we assign a {\em rank}, which is the $N+2$-tuple 
$$
\big(\sum_i e_i, j,  e_1, \ldots , e_N\big)
$$
and order the ranks lexicographically.  

Given a polynomial $f \in R\{\ybar\}$,
the {\em leader} of $f$, denoted $u_f$ is the differential variable of maximal rank that
occurs in $f$.  Writing $f$ as a polynomial in $u_f$, $f = \sum_{j=0}^d I_ju_f^j$,
the {\em initial} of $f$, denoted $I_f$ is the polynomial $I_d$.  The {\em separant}
of $f$, written $S_f$, is the derivative of $f$ with respect to $u_f$.  In other words,
$$
S_f =  (\pd / \pd u_f) f = I_1 + 2I_2u_f + \ldots + dI_du_f^{d-1}
$$
Given a finite set $A \sub R\{\ybar\}$, let $H_A$ denote the product
$\prod_{f \in A}S_fI_f \in R\{\ybar\}$.
\edf

Even when $R$ is a differential field, prime differential ideals in polynomial rings
over $R$ are
not necessarily finitely generated.  Nevertheless, we have the following fact,
which is a consequence of Ritt's division algorithm (Proposition 1 on p.\ 79
of~\cite{EK1}).

\bla
\label{lm:ideals}
Let $R\{y_1, \ldots , y_n\}$ be a differential polynomial ring, and let
$\mfp$ be a differential prime ideal with $\mfp \cap R = \{0\}$.
There is a finite set $A \sub \mfp$ (called the characteristic set of $\mfp$)
such that for any polynomial $g$, $g \in \mfp$ if and only if is an $m \in \NN$ with
$$
H_A^m\cdot g \in [A]
$$
In addition, for each $f \in A$, $I_f \not\in \mfp$ and $S_f \not\in \mfp$.
\ela

\brm
When $R$ is a field, this is a result of Rosenfeld (\cite{EK1}, p.\ 167).  Otherwise, it also
follows from Proposition 1, on p.\ 79, Lemma 8, on p.\ 82, and the Remark on p.\ 124.
\erm

\subsection{Differentially closed fields}

A differential field $K$ is called {\em differentially closed} if, for any finite set of 
differential polynomials 
$$
p_1(x_1, \ldots , x_m), \ldots, p_s(x_1, \ldots , x_m), q_1(x_1, \ldots , x_m), \ldots ,
q_t(x_1, \ldots , x_m)
$$
the following property holds.  If there is a differential field extension $L$ of $K$,
and an $m$-tuple $\abar \in L$ such that $p_i(\abar) = 0$, for all $i \leq s$, and
$q_j(\abar) \neq 0$, for all $j \leq t$, then such a tuple can be found already in $K$.
(In the language of model theory, these are precisely the existentially closed
differential fields.)  In the context of a single derivation, A. Robinson~\cite{ARob}
proved that the class of differentially closed fields can be axiomatized in 
first-order logic.  Later, L. Blum found a much simpler set of axioms that involved
only differential polynomials in one variable.  In fact, she proved a very general
result which says that, in a wide-range of contexts, if the class of existentially closed 
models of some class of algebraic structures is first-order axiomatizable, then it has a 
set of first-order axioms that only mention functions in one variable
(see~\cite{GS}, Theorem 17.2, also~\cite{CW}, Lemma 2.2).

First-order axioms for the class of differentially closed fields with some fixed finite
number of derivations were first given by McGrail~\cite{TM}.  
The next lemma is an immediate consequence of her work.

\bla
\label{lm:TM}
Let $K$ be a differential field, with $n$ commuting derivations.  Then $K$ is 
differentially closed if and only if for any prime ideal $\mfp \sub K\{x\}$ and
any $B \in K\{x\} \setminus \mfp$, there is an $a \in K$ such that $f(a) = 0$,
for all $f \in \mfp$, and $B(a) \neq 0$.
\ela

Below we will need two facts about differentially closed fields.  

\bla
Any differential field can be embedded in a differentially closed field.
\ela

\bla 
\label{lm:qe}
Let $K$ be a differentially closed field, and consider a finite set of differential 
polynomials $p_1(\xbar, y), \ldots , p_j(\xbar,y),
q(\xbar,y)$ in $\ZZ\{\xbar, y\}$, where $\xbar = (x_1, \ldots , x_n)$
is an $n$-tuple of variables.  Suppose that there is an $n$-tuple $\cbar \in K$,
and $t \in K$, such that for all $i \leq j$, $p_i(\cbar, t) = 0$, and 
$q(\cbar, t) \neq 0$.  

Then there is another set of polynomials 
$f_1(\xbar), \ldots ,f_k(\xbar), g(\xbar)$, such that:
\ben
\item for all $i \leq k$, $f_i(\cbar) = 0$ and also $g(\cbar) \neq 0$;
\item given any other $n$-tuple $\ebar \in K$, if 
for all $i \leq k$, $f_i(\ebar) = 0$ and also $g(\ebar) \neq 0$, then there is an $s \in K$,
such that for all $i \leq j$, $p_i(\ebar,s) = 0$, and $q(\ebar,s) \neq 0$.
\een
\ela

\bpr
This is a consequence of quantifier elimination~\cite{TM} 
or Seidenberg's related elimination
theorem~\cite{AS} (see also~\cite{EK3}, p.\ 578).
\epr

We now explain the geometric meaning of this lemma.  Recall that an
affine differential algebraic variety $V \sub K^m$ is the set of zeros of
some finite set of differential polynomial functions in $K\{x_1, \ldots ,x_m\}$.
These are the basic closed sets in the Kolchin topology, which plays
the role here of the Zariski topology in algebraic geometry.  
Define a {\em quasi-affine differential algebraic variety} to be an open
subset of an affine differential algebraic variety.  (In other words,
a quasi-affine variety $U \sub K^m$ is the set of points $\cbar \in V$,
such that $q(\cbar) \neq 0$, for some fixed affine variety $V$ and
some fixed polynomial $q(\xbar)$.)  The previous lemma can 
be restated in the following form.

\bla
Let $X$ be a quasi-affine 
variety, $Y$ an affine variety, and $f:X \ra Y$ a morphism.  For each
$y \in Y$, let $X_y$ denote the fiber of $X$ above $y$.
If for some generic element $y_0 \in Y$,
the fiber $X_{y_0}$ is non-empty, then there is a quasi-affine
variety $U \sub Y$ such that for all $u\in U$, the fiber $X_u$ is
non-empty.
\ela

\section{A Chevalley type theorem}

\subsection{Main result}

Our main result is an analog of a theorem of Chevalley, for rings
with commuting derivations.  The basic structure of the argument 
follows Kac's argument for the case of a single derivation,
though some of the algebra is replaced by an appeal to quantifier elimination
for differentially closed fields.  

\bthm
\label{thm:main}
Let $S$ be a differential algebra with no zero divisors, $R$ a differential
subalgebra of $S$ over which is $S$ is differentially finitely generated,
and $K$ a differentially closed field.  Then for any nonzero $b \in S$, 
there is a nonzero $a \in R$ such that any homomorphism $\phi: R \ra K$
with $\phi(a) \neq 0$ extends to a homomorphism $\psi: S \ra K$ with 
$\psi(b) \neq 0$.
\ethm

\bpr
We argue by induction on the minimal number $n$ of elements that generate
$S$ over $R$.  Suppose first that we have proved the base case $n = 1$,
and consider $S = R\{t_1, \ldots , t_n\}$, $n > 1$, and nonzero $b \in S$.
Let $R_1 = R\{t_1, \ldots , t_{n-1}\}$.  By hypothesis, there is an $a_1 \in R_1$
such that any homomorphism $\phi_1: R_1 \ra K$ with $\phi_1(a_1) \ne 0$
lifts to a homomorphism $\psi: S \ra K$ with $\psi(b) \ne 0$.  Likewise, there
is an $a \in R$ such that any homomorphism $\phi: R \ra K$ with $\phi(a) \ne 0$,
lifts to a homomorphism $\psi_1: R_1 \ra K$ with $\psi_1(a_1)\ne 0$.
Together these two claims imply that $a \in R$ is as desired.  
It only remains then to establish the case $n=1$.

We now assume that $S = R\{t\}$ and $b = B(t)$, for some nonzero polynomial
$B(y) \in R\{y\}$.  Given a homomomorphism $\phi: R \ra K$ and a polynomial
$f(y) \in R\{y\}$, let $f^\phi(y) \in K\{y\}$ be the polynomial obtained by applying
$\phi$ to each coefficient in $f(y)$.
There are two cases to consider.

\msk

\noi
\underline{Case 1.}  The element $t$ is differentially transcendental over $R$.
Choose $a \in R$ to be any nonzero coefficient of $B(y)$.  Let $\phi:R \ra K$
be any homomorphism with $\phi(a) \ne 0$.  As $t$ is differentially 
transcendental, for any $c \in K$, there is a unique homomorphism 
from $S$ to $K$ lifting $\phi$ that sends $t$ to $c$.  Since $K$ is 
differentially closed, there exists a $c_0 \in K$ with $B^\phi(c_0) \neq 0$.
We can then choose $\psi: S \ra K$ to be the homomorphism lifting 
$\phi$ such that $\psi(t) = c_0$.

\msk

\noi
\underline{Case 2.}
Finally, suppose $t$ is differentially algebraic over $R$.
Let $\mfp \sub R\{y\}$ be the differential prime ideal $\mfp = \{f(y): f(t) = 0\}$.
By Lemma~\ref{lm:ideals}, there is a finite set $A \sub \mfp$, such that for any polynomial $f$,
$f \in \mfp$ if and only if  there is an $m$, such that $H_A^m \cdot f \in [A]$.
Write $A = \{p_1(y), \ldots , p_j(y)\}$.

Let $\cbar$ be an enumeration of all the coefficients that occur in the polynomials
$B(y), H_A(y), p_1(y), \ldots , p_j(y)$, and let $M$ be the size of $\cbar$.  We will write
$\hat{B}(\cbar,y)$, $\hat{H}_A(\cbar,y)$, and $\hat{p}_i(\cbar,y), i \leq j$, in order to exhibit 
all of the coefficients explicitly.  Then
$\hat{B}(\zbar,y)$, $\hat{H}_A(\zbar,y)$, and $\hat{p}_i(\zbar,y)$ are differential polynomials
over $\ZZ$ in $M+1$ variables.

Let $L$ be a differentially closed field containing $S$.
By Lemma~\ref{lm:qe}, there are polynomials $f_1(\zbar), \ldots, f_k(\zbar), g(\zbar)$ in 
$\ZZ\{\zbar\}$, in $M$ variables, such that:
\ben
\item  for all $i \leq k$, $f_i(\cbar) =0$ and $g(\cbar)\neq 0$;
\item  given any $M$-tuple $\ebar \in L$,
if for all $i \leq k$, $f_i(\ebar) =0$ and $g(\ebar)\neq 0$, then there is an
element $s \in L$ such that for all $i \leq j$, $\hat{p}_i(\ebar,s) = 0$,
and both $\hat{B}(\ebar,s)\neq 0$ and $\hat{H}_A(\ebar,s)\neq 0$.
\een

We can now define the desired $a \in R$ to be $g(\cbar) \in R$.  It remains to 
show that $a$ has the desired property.

Let $\phi$ be a homomorphism from $R$ to $K$.  Since for all $i$,
$f_i(\cbar) = 0$, then $f_i(\phi(\cbar)) = 0$,
where $\phi(\cbar)$ is the tuple obtained by applying $\phi$ to each
$c$ in $\cbar$.  Since $g(\phi(\cbar)) = \phi(g(\cbar))$, it is also clear that
$g(\phi(\cbar)) \neq 0$.  Again by the definition of the polynomials $f_i$ and $g$,
there is an element $v \in K$ such that  
$\hat{B}(\phi(\cbar),v) \neq 0, \hat{H}_A(\phi(\cbar),v) \neq 0$
and for all $i , \hat{p}_i(\phi(\cbar),v) = 0$.  In particular, this implies that for any
 polynomial $s(y) \in [A]$, $s^\phi(v) = 0$.  

Finally, I claim that there is a unique homomorphism $\psi: S \ra K$,
lifting $\phi$, with $\psi(t) = v$.  In order to establish this, it suffices to 
show that for every $h(y) \in \mfp$, $h^\phi(v) =0$.  
By Lemma~\ref{lm:ideals}, there is an $m \in \NN$ such that
$H_A^m(y)h(y) \in [A]$.  Write $s(y) = H_A^m(y)h(y)$.  By the previous remark,
$s^\phi(v) = 0$.  Since $s^\phi(v) = (H_A^m)^\phi(v)h^\phi(v)$ and
$H_A^\phi(v) \ne 0$, we get that $h^\phi(v)=0$, as desired.
\epr

\bco
\label{cor:char}
Let $K$ be a differentially closed field and $S$ be a finitely generated differential
$K$-algebra with no zero divisors.  For all nonzero $b \in S$, there is a homomorphism
$\phi:  S \ra K$ with $\phi(b) \ne 0$.
\eco

As a corollary, we also get a new proof of the following theorem of Kolchin
(\cite{EK1}, p.\ 140; see also  \cite{EK2}, p.\ 579), which generalized earlier
results of Ritt~\cite{JFR}, Seidenberg~\cite{AS}, and Rosenfeld~\cite{ARos}.

\bco
[Kolchin]  
\label{cor:kol}
Let $S$ be a differential ring with no zero divisors and $R$ a subring 
over which $S$ is differentially finitely generated.  Given a nonzero $b \in S$,
there is an $a \in R$ such that for every differential prime ideal $\mc{P} \sub R$
with $a \not\in \mc{P}$, there is a differential prime ideal $\mc{Q} \sub S$, not containing 
$b$, with $\mc{Q} \cap R = \mc{P}$.
\eco

\bpr
Choose $a \in R$ to be the element whose existence is guaranteed by 
Theorem~\ref{thm:main}.  Let $\mc{P}$ be a prime differential ideal in $R$,
$a \not\in \mc{P}$, so $R/\mc{P}$ is a subring of some differentially closed field $K$.
Let $\phi: R \ra K$ be the canonical map, so $\phi(a) \ne 0$. 
Then there is a homomorphism $\psi: S \ra K$, with $\psi(b) \ne 0$.
Choose $\mc{Q} = \ker \psi$.
\epr

\subsection{Ritt schemes and geometric Chevalley}

In this section, we prove a generalization of Buium's geometric differential 
Chevalley theorem~\cite{Bui1}, extending his result for a single derivation
to the case of commuting derivations.  We begin by recalling the
notion of a Ritt scheme, which is a differential analog of the usual notion
of a scheme from algebraic geometry.  With this background in place,
our result then follows rather easily from Theorem~\ref{thm:main} above.
(As an aside, it is perhaps worth noting that there are a few proposed
notions of a `differential scheme' in the literature.  See also, for example,
Buium's paper~\cite{Bui2}, where he introduces a different idea, and recent work
of J. Kovacic~\cite{Kov}.)

\bdf
Let $R$ be a Ritt algebra.  Let $\Specd R$ be the topological space
whose points are the prime differential ideals of $R$, with the topology
induced by the inclusion map $i: \Specd R \ra \Spec R$.  

One defines a sheaf of rings $\mcO$ on $\Specd R$ exactly as in the
classical case.  Then for each open $U \sub \Specd R$, $\mcO(U)$ is 
naturally a 
differential ring and, for each differential prime ideal $\mfp \sub R$,
the stalk $\mcO_\mfp$ is a local differential ring.

Call $\Specd R$, endowed with this sheaf of differential rings, an
{\em affine Ritt scheme}.
\edf

As in the classical case, given any differential ideal $\mfa \sub R$,
let $V(\mfa)$ denote $\{\mfp: \mfa \sub \mfp)$, a closed subset of $\Specd R$.

\bdf
Let $X,Y$ be affine Ritt schemes.  Say that a map $f: X \ra Y$
is of differential finite type if there are Ritt algebras $R,S$,
with $X \cong \Specd S$, $Y \cong \Specd R$, such that $f$ is induced
by a homomorphism $\phi: R \ra S$ and $S$ is differentially
finitely generated over $\phi(R)$.
\edf

We also need the following theorem of Kolchin (\cite{EK1}, p.\ 126),
which we restate in a form convenient for our application here.
It is a general form of the Ritt-Raudenbush basis theorem, the
analog for differential rings of Hilbert's basis theorem.

\bthm [Kolchin]
Let $S$ be a differential ring, $R$ a differential subring with $S$
differentially finitely generated over $S$.  Suppose that
$\Specd R$ is noetherian.  Then so is $\Specd S$.
\ethm

Recall that, given a subset of a topological space is said to be
locally closed if it is the intersection of a closed set with an open set,
and is said to be constructible if it is a finite union of locally closed sets.
We can now state our generalization of Buium's theorem.

\bthm
\label{thm:geo}
Let $f: X \ra Y$ be a morphism of differential finite type between
affine Ritt schemes.  Suppose that $Y$ is noetherian.  Then 
$f(X)$ is a constructible subset of $Y$.  
\ethm

\bpr
Let $X \cong \Specd S$, $Y \cong \Specd R$, and let $f$ be
the map induced by the homomorphism $\phi: R \ra S$.
To prove that $f(X)$ is constructible, it suffices to show that for 
any irreducible closed set $Y' \sub Y$, $f(X) \cap Y'$
either contains a non-empty open subset of $Y_0$ or  is
not dense in $Y'$ (e.g.,~\cite{Mat}, Proposition (6.C)).
Let $Y'$ be such a set, so $Y' = V(\mfq)$, for some prime
differential ideal $\mfq \sub R$.  Define $X' = f^{-1}(Y')$,
and let $f'$ be the restriction of $f$ to $X'$. 
Then $f': X' \ra Y'$ is the morphism induced by 
the canonical homomorphism $\phi': R/ \mfp \ra S / S\mfp$.
Thus we have reduced to the case where $R$ has
no zero divisors.

Since $X$ is noetherian, so is $X'$, which implies that 
$S/ S\mfp$ has finitely many minimal prime differential ideals,
$\mfq_1, \ldots , \mfq_n$, which correspond to the maximal
closed irrreducible components $X_1, \ldots , X_n$ of $X'$,
$X_i = V(\mfq_i)$.  For all $i$, let $f_i$ be the restriction of $f'$
to $X_i$, so that $f_i: X_i \ra Y$ is induced by the canonical 
homomorphism $\phi_i: R/ \mfp \ra S / \mfq_i$.  
Thus we may assume that both $R$ and $S$ have no zero divisors.

It now suffices to show, for each $i$, that either $f_i(X_i)$ is not dense
in $Y$, or it contains a nonempty open subset of $Y$.
Suppose first that $f_i$ is not injective, so $\ker (f_i) = \mfp'$, with
$\mfp'$ a non-trivial prime differential ideal of $R$.  Then 
$f_i(X_i) \sub V(\mfp')$, a proper closed subset of $Y$.
Otherwise, suppose that $f_i$ is injective.
By Corollary~\ref{cor:kol}, there is a $b \in R$ such that
for any differential prime ideal $\mfp \sub R$, with $b \not\in \mfp$,
there is a differential prime ideal $\mfq \sub S$, with 
$\mfq \cap R = \mfp$.  Equivalently $f(\mfq) = \mfp$.  
In particular, this shows that $f(X)$ contains the basic open
set $D(b) = \Specd R \setminus V(\{b\})$.
\epr

\subsection{Differential fields with quantifier elimination}

In this section, we briefly discuss the relationship between 
(various versions of) the differential Chevalley theorem and 
quantifier elimination for differentially closed fields.  Expressed somewhat
informally, Theorem~\ref{thm:geo} says that if $f:X \ra Y$ is a morphism of
(affine) differential varieties, then $f(X)$ is a constructible set.  
Quantifier elimination for differentially closed fields 
can be formulated as the statement that, under the same hypotheses,
$f(X(K))$ --- the image of the $K$-valued points of $X$ under $f$ --- is 
a constructible subset of $Y(K)$.  Despite the similarities
between these statements, neither one can be deduced
directly from the other.  In particular, the geometric version of the Chevalley
theorem holds for varieties over arbitrary differential fields,
so it is certainly more general than quantifier elimination.  In the other
direction, geometric Chevalley does not immediately imply quantifier 
elimination, since one cannot, e.g.,  replace $K$ by any other differential field.

Nevertheless, quantifier elimination does follow from the more
algebraic Theorem~\ref{thm:main}, since morphisms 
$\phi: R \ra K$ correspond precisely to $K$-valued points of the 
variety $\Specd R$.  One can then establish quantifier elimination by
arguing as in the proof of Theorem~\ref{thm:geo}.  On the other hand,
there is no reverse implication as,
\ben
\item [(i)]  by Theorem~\ref{thm:3} below, the conditions for 
Theorem~\ref{thm:main} yield a characterization of differentially closed fields;
\item [(ii)]  it is known that (with a single derivation)  there are differential fields 
with quantifier elimination that are not differentially closed~\cite{HI}.
\een
The following related question remains open.

\bqu
Given a fixed number $n > 1$ of commuting derivations, are there 
differential fields with quantifier elimination that are not differentially
closed?
\equ

\section{Characterizing differential fields}

In this section, we provide a number of equivalent characterizations of
differentially closed fields.  For the most part, we follow the proof of 
Theorem 3 of~\cite{Kac}.  We should also mention that Kolchin's
notion of a constrainedly closed field provides another description
(see~\cite{EKcon}).

\bpro
\label{pr:dcf.eq}
Let $K$ be a differential field.  The following properties are equivalent.

\ben
\item  $K$ is differentially closed.
\item Let $\mfp$ be a prime differential ideal of $K\{y\}$, and $B \in K\{y\} \setminus \mfp$.
Then there exists $a \in K$ such that $f(a) = 0$, for all $f \in \mfp$, and $B(a) \ne 0$.
\item  Let $\mfp$ be a prime differential ideal of $K\{y_1, \ldots, y_m\}$, and 
$B \in K\{y_1, \ldots , y_m\} \setminus \mfp$.  Then there exists an $m$-tuple $\abar$ in $K$ such that 
$f(\abar) = 0$, for all $f \in I$, and $B(\abar) \ne 0$.
\item  Let $\mfp$ be a prime differential ideal of $K\{y_1, \ldots , y_m\}$, $m \geq 1$.
Then there exists an $m$-tuple $\abar$ in $K$ such that $f(\abar) = 0$, for all $f \in \mfp$.
\item  Let $\mfm$ be a maximal differential ideal of $K\{y_1, \ldots , y_m\}$, $m \geq 1$.
Then there exists an $m$-tuple $\abar$ in $K$ such that $f(\abar) = 0$, for all $f \in \mfm$.
\een
\epro

\bpr   
By Lemma~\ref{lm:TM}, (1) and (2) are equivalent.
We establish that (1) implies (4).  Given $\mfp$, let $F$ be the differential
field $K\{y_1, \ldots , y_m\} / \mfp$ and let $b_i = y_i + \mfp \in F$, for all $i \leq m$.
Let $A$ be the characteristic set of $\mfp$, $A = \{f_1(\ybar), \ldots , f_k(\ybar)\}$.
Since $f_i(\bbar) = 0$, for all $i$, and $H_A(\bbar) \neq 0$, there is a tuple 
$\cbar$ in $K$ such that $f_i(\cbar) = 0$, for all $i$, and $H_A(\cbar) \neq 0$.
It remains to show that for all $h(\ybar) \in \mfp$, $h(\cbar) = 0$.  
By Lemma~\ref{lm:ideals}, there is an $n$ such that $H_A^n(\ybar)h(\ybar) \in [A]$, which
implies that $H_A^n(\cbar)h(\cbar) = 0$ and, thus, $h(\cbar) = 0$.

To show (4) implies (3), assume that $\mfp$ is a prime differential ideal of $K\{y_1, \ldots, y_m\}$, 
and $B \in K\{y_1, \ldots , y_m\} \setminus \mfp$.  As above let  $F = K\{y_1, \ldots , y_m\} / \mfp$
and let $b_i = y_i + \mfp \in F$, for all $i \leq m$.  Since $B(\bbar) \ne 0$, we can set
$b' = 1 / B(\bbar)$.  Let $\mfq \sub K\{y_1, \ldots, y_m, y_{m+1}\}$ be the prime differential
ideal consisting of all polynomials $g$ such that $g(\bbar,b')=0$.  In particular,
$B(\ybar)y_{m+1} - 1 \in \mfq$.  By hypothesis, there is a tuple $(\cbar,c')$ in $K$,
with $g(\cbar,c') = 0$, for all $g \in \mfq$.  It is easy to see that $\cbar$ is the desired $m$-tuple.

Trivially (3) implies (2).  Finally, since every prime differential ideal embeds in a maximal
differential ideal, which is necessarily prime, we get the equivalence of (4) and (5).
\epr

The next theorem demonstrates that the Theorem~\ref{thm:main} and 
Corollary~\ref{cor:char} yield characterizations of differentially closed fields.
In the context of a single derivation, Kac proved the equivalence of (1) and (3).
The equivalence with (2) is new here.

\bthm
\label{thm:3}  Let $K$ be a differential field.  The following properties are
equivalent.
\ben
\item  $K$ is differentially closed.
\item  Let $S$ be a differential ring with no zero divisors, $b$ a nonzero element
of $S$, and $R$ a subring of $S$ over which $S$ is finitely generated.
Then there is an $a \in R$ such that for any homomorphism $\phi: R \ra K$,
with $\phi(a) \ne 0$, there is a homomorphism $\psi:S \ra K$ extending 
$\phi$ such that $\psi(b) \neq 0$.
\item  For any finitely generated differential $K$-algebra $S$ with no zero divisors and 
nonzero $b \in S$, there is a
$K$-algebra homomorphism $\psi: S \ra K$ with $\psi(b) \ne 0$.
\item
For any finitely generated differential $K$-algebra $S$ with no zero divisors
there is a $K$-algebra homomorphism $\psi: S \ra K$.
\item
For any finitely generated differential $K$-algebra $S$
there is a $K$-algebra homomorphism $\psi: S \ra K$.
\een
\ethm

\bpr
By Theorem~\ref{thm:main}, (1) implies (2).  Letting $R = K$, one gets that
(2) implies (3).

To prove (3) implies (1) suppose that $K$ is not differentially closed.  
Then there are polynomials
$p_1(\xbar), \ldots, p_m(\xbar), q(\xbar)$ in $K\{\xbar\}$ with the following
property. There is a tuple $\sbar$ in some differential field extension $L$ of $K$,
with $p_i(\sbar) = 0$, for all $i$, and $q(\sbar) \ne 0$, but there is no such 
tuple in $K$ itself.  Let $S$ be the $K$-algebra $K\{\sbar\} \sub L$
and  let $b = q(\sbar)$.  
Suppose now that $\phi: S \ra K$ is a $K$-algebra homomorphism,
and let $\sbar ' = \phi(\sbar)$.  Clearly $p_i(\sbar ') = 0$, for all $i$,
so by assumption $q(\sbar ') = 0$.  Since $\phi(b) = q(\sbar ')$, we are done.

Immediately, (3) implies (4).  To prove the reverse implication, given
$S$ and $b$, apply (4) to the ring $S[b^{-1}]$.  It is also immediate that
(5) implies (4).  In the other direction, given a differential $K$-algebra $S$,
let $\mfp$ be a minimal prime ideal, which is necessarily  a 
differential ideal~\cite{Gil}.  Let $S' = S / \mfp$ and let 
$\theta: S \ra S'$ be the canonical map.  By assumption, there is a homomorphism
$\eta: S' \ra K$.  We can then define $\psi = \eta \circ \theta$.
\epr

\brm
\label{rm:main}
\ben
\item
In the statement of the previous theorem, one may replace (2) by property (2$'$),
which requires additionally that both $S$ and $R$ be finitely generated.  Indeed, 
to show that (2$'$) implies (1), proceding as in the above proof, let $R$ be the 
differential $\QQ$-algebra generated by the set of all coefficients that occur in the 
polynomials  $p_1, \ldots , p_m, q$, and let $S = R\{\sbar\}$.  Then by the 
previous argument, for any homomorphism $\psi : S \ra K$ that lifts the map
$\phi$ embedding $R$ into $K$, $\psi(b) = 0$.

\item
Further, by the axioms for differentially closed fields, or Proposition~\ref{pr:dcf.eq}, 
one may also assume in (2), resp.\ (3), that $S$ is generated by a single element 
over $K$, resp.\ $R$.
\een
\erm

The theorem can be stated more geometrically.

\bco  Let $K$ be a differential field.  The following properties are equivalent.
\ben
\item $K$ is differentially closed.

\item For any (affine) differential
algebraic variety $V$ over $K$ and any nonzero regular function $f$ on $V$, there 
is a $v \in V(K)$ with $f(v) \ne 0$.
(In fact, one may restrict attention to $V \sub \mathbb{A}^1$.)

\item Any (affine) differential
algebraic variety $V$ over $K$ has a $K$-rational point.
\een
\eco

\bco[Differerential Nullstellensatz]
\label{diff.null}
Let $K$ be a differentially closed field, an let $\mc{I}$ be a proper differential ideal in 
$K\{x_1, \ldots, x_n\}$.  Then $\mc{I}$ has a zero in $K$.
\eco

\bpr
Let $S = K\{x_1, \ldots, x_n\} / \mc{I}$ and, by (5) of Theorem~\ref{thm:3}, let
$\psi: S \ra K$ be a homomorphism.  
Then $\psi(\xbar)$ is a zero of $\mc{I}$.
\epr

\bex
We describe a simple example of a differential variety $V$ over a differential field $K$
with a nonzero regular function on it which is nevertheless identically zero on 
all points in $V(K)$.  Without loss of generality, we assume that $K$ has a single derivation.
Recall that the set of solutions of a homogeneous linear differential equation of order $n$,
$\pd^{(n)}y + a_{n-1}\pd^{(n-1)}y + \ldots + a_1 \pd y + a_0 = 0$, over some differential 
field $K$ is a vector space of dimension $\leq n$ over the field $k \sub K$ of constants,
$k = \{a \in K : \pd (a) = 0\}$.  Moreover, there is always a field extensions $L$ of $K$ in 
which the solution set has dimension $= n$.  (For example, see~\cite{PS}.)  

Choose such an equation and a differential field $K$
such that the dimension of the space of solutions, $W \sub K$, is $< n$.  
Let $V$ be the differential variety $V = W^n \sub K^n$.  Let $Wr: K^n \ra K$ be the Wronskian, 
which is a differential polynomial function with the property that $Wr(a_1, \ldots, a_n) = 0$ if 
and only if the $a_i$ are linearly dependent over the constants.  By the choice of $K$, 
$Wr$ is identically zero on $V$, but this does not remain true when passing to an
extension field $L$ in which the solution set of the original equation has maximal dimension.
 \eex

\section{Large differentially closed fields}
\label{sec:sat}

In the statement of Theorem~\ref{thm:main}, the differential algebra $S$ is required to be 
finitely generated over $R$.  An easy example, below, shows that this condition
is necessary.  We then prove a version of the theorem that loosens this
restriction.

\bex
Let $R = \ZZ\{x_i : i \in \NN\}$ and let $S= R\{x_i^{-1}: i \in \NN\}$.  
Given any differentially closed field $K$ and any $b \in R$, 
let $\phi: R \ra K$ be a homomorphism with $\phi(b) \ne 0$
and $\phi(x_i) = 0$, for some $x_i$ not occuring in $b$.
Since $x_i$ is invertible in $S$, $\phi$ cannot be extended to 
a homomorphism $\psi: S \ra K$.
\eex

To begin, we prove a differential version of the lying over and going up theorem
for integral ring extensions, which may be of independent interest.  
Recall that, given a ring $S$, a subring $R \sub S$, and ideals $I \sub R$,
$J \sub S$, then $J$ is said to be lying over $I$ if $J \cap R = I$.
Our proof uses the following result (\cite{IK}, p.\ 23).

\bpro
Let $S$ be a differential ring, $R$ a subring of $S$, and $\mfp \sub R$ 
a prime ideal.  Suppose that $I$ is a radical differential ideal lying over $\mfp$
and, for all $r \in R, s \in S$, if $rs \in I$, then $r \in I$ or $s \in I$.  
Then $I$ is the intersection of prime differential ideals lying over $\mfp$.
\epro

\bthm[Differential lying over and going up theorem]
\label{thm:up}
Let $S$ be a differential ring and $R$ a subring, with $S$ integral over $R$.
Then for any differential prime ideal $\mfp \sub R$, there is a differential 
prime ideal $\mfq \sub S$ such that $\mfq \cap R = \mfp$.  Further, $\mfq$
can be chosen so as to contain any differential ideal $\mfq_0$, with 
$\mfq_0 \cap R \sub \mfp$.

\ethm

\bpr
We begin by localizing at $\mfp$, so let $R' = R_\mfp$,
$S' = S_\mfp = S \otimes_R R_\mfp$, and $\mfp ' = \mfp R_\mfp$.
One can easily check that $S'\mfp'$ is a differential ideal, so by a 
remark above, its radical $\{S'\mfp'\}$ is also a differential ideal, which also 
lies above $\mfp'$.  In order to apply the previous proposition, we need to 
verify that, given $r \in R',s \in S'$, if $rs \in \{S'\mfp'\}$, then $r \in \{S'\mfp'\}$
or $s \in \{S'\mfp'\}$.  But this follows from the fact that $r$ is either a unit or
in $\mfp'$, because $R'$ is a local ring.  Thus $\{S'\mfp'\}$ is the
intersection of differential prime ideals lying above $\mfp'$.
Let $\mfq'$ be one such ideal, and let $\mfq$ be its preimage under
the canonical map $S \ra S'$.  Then $\mfq$ is as desired.
To prove the second statement, one first replaces $S$ with $S/\mfq_0$,
and then argues as before.
\epr

Given an infinite cardinal $\lambda$, there is a model theoretic notion of a 
field (or any algebraic structure) being $\lambda$-{\em saturated}.  When the 
field has quantifier elimination, one has an easy to state equivalent
algebraic condition, which involves the existence of solutions
to sets of polynomial equalities and inequalities of bounded size.

\bla  Let $\lambda$ be an infinite cardinal.  
\ben  
\item  If $\lambda = \aleph_0$, the countably infinite cardinal, a 
differentially closed field is $\lambda$-saturated if and only if given any
set of polynomials
$$
\{p_i(x): i \in I\} \cup \{q_j(x) : j \in J\} \sub k\{x\}
$$
with $k$ a finitely generated differential subfield of $K$,
the following property holds.
If there is a differential field extension $L$ of $K$ and an $a \in L$ such that 
$p_i(a) = 0$, for all $i \in I$, and $q_j(a) \neq 0$, for all $j \in J$,
then there is such an $a$ already in $K$.
\item
If $\lambda$ is uncountable, a differentially closed field is $\lambda$-saturated
if and only if given any set of polynomials 
$$
\{p_i(x): i \in I\} \cup \{q_j(x) : j \in J\} \sub k\{x\}
$$
with $k$ a differential subfield of $K$ of cardinality less than $\lambda$, 
the following property holds.
If there is a differential field extension $L$ of $K$ and an $a \in L$ such that 
$p_i(a) = 0$, for all $i \in I$, and $q_j(a) \neq 0$, for all $j \in J$,
then there is such an $a$ already in $K$.
\een
\ela

\bpr
This follows immediately from quantifier elimination.
\epr

\brm
It is well-known that in the statement of the previous lemma, one may also allow sets of 
polynomials over any fixed number of variables.
\erm

One can show that for all $\lambda$, there is a $\lambda$-saturated 
differentially closed field of cardinality $\kappa$, for every cardinal $\kappa \geq \lambda$.  
This can be proved using Zorn's Lemma or, alternatively, using some model theory, it is a 
consequence of the fact that the theory of differentially closed fields is $\omega$-stable.
We will also use the following lemma, which gives a somewhat different characterization of
saturated differentially closed fields.

\bla
\label{lm:sat}  Let $\lambda$ be an infinite cardinal and $K$ a differentially closed field.  
The following are equivalent.
\ben
\item  $K$ is $\lambda$-saturated.

\item  Let $L$ a differential field of cardinality $\leq \lambda$ and $F$ a subfield of $L$ of 
cardinality $<\lambda$, if $\lambda$ is uncountable, and finitely generated otherwise.
Then any embedding of $F$ into $K$ extends to an embedding of $L$ into $K$.
\een
\ela

\bpr
This is a consequence of a standard model theoretic fact, but 
one can also give an easy direct proof.
\epr

\bthm
\label{thm:card} 
Let $K$ be a $\lambda$-saturated differentially closed field, $\lambda$ an infinite 
cardinal.  Let $S$ be a differential ring of cardinality $\leq \lambda$ with no zero divisors, 
$b$ a nonzero element of $S$,
and $R$ a subring of cardinality $< \lambda$, if $\lambda$ is uncountable, and
finitely generated otherwise.  Suppose that $S$ is integral over some differential
subring $T$ which is finitely generated over $R$.  Then there is an $a \in R$
such that any homomorphism $\phi: R \ra K$ with $\phi(a) \ne 0$ extends to a
homomorphism $\psi: S \ra K$ with $\psi(b) \ne 0$.
\ethm

\bpr
Let $S' = T\cup\{b\}$.  By Theorem~\ref{thm:main}, choose $a \in R$
such that any homomorphism from $R$ to $K$ that does not annihilate $a$
extends to a homomorphism from $S'$ to $K$ that does not annihilate $b$.
Suppose now that we have $\phi: R \ra K$ with $\phi(a) \ne 0$.  
By hypothesis, there is a map $\theta: S' \ra K$ with $\theta(b) \ne 0$.  
Let $\mfp = \ker \theta$.  Since $\mfp$ is a differential prime ideal,
by Theorem~\ref{thm:up} there is a differential prime ideal $\mfq \sub S$
lying over $\mfp$.  Let $F$ be the fraction field of $S'/\mfp$ and $K$ be the
fraction field of $S/\mfq$, which extends $F$.  The map $\theta$ extends
uniquely to a homomorphism $\theta_0: F \ra K$ so, by Lemma~\ref{lm:sat},
there is an embedding $\psi_0:L \ra K$ extending $\theta_0$.  We can then
choose $\psi$ to be the composition of the canonical map $S \ra S/ \mfq$
with $\psi_0$. 
\epr

We also get the following characterization of saturated differentially closed fields
in analogy to Theorem 3 of~\cite{Kac}.

\bthm  Let $\lambda$ be an infinite cardinal and let $K$ be a differential field.  
The following are equivalent.
\ben
\item  $K$ is a $\lambda$-saturated differentially closed field.
\item  
Let $S$ be a  differential algebra of cardinality $ \leq \lambda$ 
with no zero divisors and let $R$ be a differential 
subalgebra of $S$ of cardinality $< \lambda$, if $\lambda$ is uncountable, and
finitely generated otherwise.  
Suppose also that $S$ is integral over $R$.  Let $b$ be a nonzero element of $S$.  
Then there is an $a \in R$ such that any homomorphism $\phi: R \ra K$,
with $\phi(a) \neq 0$, lifts to a homomorphism $\psi:S \ra K$, with $\psi(b) \neq 0$.
\een
\ethm

\bpr
That $(1)$ implies $(2)$ is the content of Theorem~\ref{thm:card}.

In the other direction, let us assume that $\lambda$ is uncountable.
(The argument for $\lambda$ countable is similar.)
If $K$ is not a differentially closed field, then by Remark~\ref{rm:main}, 
(2) already fails for some countable $S$ and $R$, both finitely generated.   So suppose that 
$K$ is a differentially closed field that is not $\lambda$-saturated.  
In particular, there is a set of polynomials
$$
\{p_i(x): i \in I\} \cup \{q_j(x) : j \in J\} \sub k\{x\}
$$
with $k$ a differential subfield of $K$ of cardinality $< \lambda$, with the following property.
There is a differential field extension $L$ of $K$ and a $b \in L$ such that 
$p_i(b) = 0$, for all $i \in I$, and $q_j(b) \neq 0$, for all $j \in J$, but there
is no such $b$ in $K$ itself.

Let $R$ be the fraction field of the differential subring of $K$ generated by the coefficients 
in all the $p_i(x)$ and $q_j(x)$, and let $S$ be the fraction field of the differential subring of $L$ 
generated by  $R \cup \{b\}$.  Let $\phi$ be the inclusion map from $R$ to $K$,
so that for all nonzero $a \in R$, $\phi(a) \neq 0$.  We claim that there is no 
extension at all of $\phi$ to a map $\psi: S \ra K$.  Indeed, suppose for contradiction that
such a $\psi$ existed.  Since $S$ is a field, by the choice of $b$ this would imply 
that for all $i \in I$, $p_i(\psi(b)) = 0$, and for all $j \in J$, $q_{j}(\psi(b)) = 0$.
But this contradicts our assumption that there is no element such as $\psi(b)$ in $K$.
\epr

\section{Differential fields in positive characteristic}
\label{sec:char.p}

In this section, we consider differential fields in characteristic $p$
with $N$ commuting derivations, for fixed $p$ and $N$.  For such a field $K$,
any derivation $\dl$ is trivial on $K^p$, the field of $p^{th}$ powers, since 
$\dl(a^p) = pa^{p-1}\dl(a) = 0$, for any $a \in K$.  In particular, if $K$ is 
algebraically closed, then any derivation is trivial.  On the other hand,
given a separable extension $L$ of $K$, a derivation on $K$ can be extended in a 
unique way to a derivation on $L$, and commuting derivations extend to commuting
derivations~(\cite{EK1}, p.\ 90).  This implies that any differentially closed field is
separably closed.  
For $N = 1$, the model theory of such fields has been analyzed by Wood~\cite{Wood},
though we will not use any of her results directly.  Bu it is worth mentioning
the fact that the first-order theory of such fields does not have quantifier elimination,
as this observation is implicit in our proof that the analog of Theorem~\ref{thm:main}
fails in this context.  For $N > 1$, there is recent work of Pierce~\cite{Pierce}.

\bpro
\label{pr:dcfp}
For any prime $p$, there is a finitely generated differential ring $S$ of characteristic $p$
with no zero divisors, a differential subring $R \sub S$, and a differentially closed
field $K$ of characteristic $p$, with the following property.  There is an 
embedding $\phi: R \ra K$ that cannot be extended to a homomorphism
$\psi: S \ra K$.
\epro

In particular, Theorem~\ref{thm:main} fails in positive characteristic.

\bpr
Let $S = \FF_p(t)$, equipped with the trivial derivation $\dl$, and let 
$R = \FF_p(t^p)$.  Let $K$ be any differentially closed field of characteristic $p$
that contains a differentially transcendental element $x \in K$, such that $\dl(x) \ne 0$.
Let $\phi: R \ra K$ be the unique homomorphism with $\phi(t^p) = x^p$.
The only {\em ring} homomorphism $\psi: S \ra K$ that extends $\phi$
maps $t$ to $x$.  But $\dl(\psi(t)) \ne \psi(\dl(t)) = 0$, so $\psi$ is not a 
differential ring homomorphism.
\epr

We now prove a restricted form of the Theorem~\ref{thm:main}
in positive characteristic, where the subring $R$ is identified with $K$.  

\bthm
\label{thm:main.p}
Let $K$ be a differentially closed field of characteristic $p$ with $N$ commuting
derivations.  Let $S$ be a finitely generated $K$-algebra with no zero divisors, 
and let $b$ be a  nonzero element in $S$.  Then there is a $K$-algebra homomorphism
$\psi: S \ra K$ with $\psi(b) \neq 0$.
\ethm

\bpr
Let $\cbar$ be a finite tuple of elements in $S$ that generates $S$ over $R$,
and let $b = g(\cbar)$, for some polynomial $g(\xbar) \in K\{\xbar\}$.
Define $\mfp \sub K\{\xbar\}$ to be set of $f(\xbar)$ such that $f(\cbar)=0$,
which is a prime differential ideal.  Suppose first that $\mfp = 0$.  In this case,
it suffices to find a tuple $\ebar \in K$ such that $g(\ebar) \ne 0$, which is 
possible since $K$ is differentially closed.  We can then define 
$\psi$ as the unique homomorphism sending $\cbar$ to $\ebar$.  

Otherwise, by Lemma~\ref{lm:ideals}, there is a finite set $A \sub \mfp$ such that for any
polynomial $g \in \mfp$ if and only if there is an $m \in \NN$ with 
$H_A^m\cdot g \in [A]$.  Furthermore, $H_A(\cbar) \neq 0$.  
Let $L$ be the fraction field of $S$, which can be 
equipped in a unique way with commuting derivations agreeing with those on $S$.
Thus, the set of polynomial equalities and inequalities,
$$
\{f(\xbar) = 0 : f \in A\} \cup \{H_A(\xbar) \neq 0, g(\xbar) \neq 0\} 
$$
has a solution in a field extending $K$.
Since $K$ is differentially closed, there is a tuple $\ebar \in K$ 
satisfying the same set of equalities and inequalities.  

It remains to show that the mapping $\cbar \mapsto \ebar$ extends to a homomorphism
$\psi: S \ra K$.   For then $\psi(b) = \psi(g(\cbar)) = g(\ebar) \neq 0$.
To do so, it suffices to prove that for all $q \in \mfp$, $q(\ebar) = 0$.
By Lemma~\ref{lm:ideals}, again, there is an $m$ such that $H^m_A\cdot q \in [A]$,
so $H^m_A(\ebar)q(\ebar) = 0$.  By the choice of $\ebar$, $H^m_A(\ebar) \neq 0$,
so $q(\ebar) = 0$, as desired.
\epr

As in characteristic 0, we obtain a new characterization of differentially closed fields.

\bthm  Let $K$ be a differential field of characteristic $p$.  The following properties
are equivalent.
\ben
\item  $K$ is differentially closed.

\item  For any finitely generated differential $K$-algebra $S$ with no zero divisors
and any nonzero $b \in S$, there is an algebra homomorphism $\psi: S \ra K$
with $\psi(b) \neq 0$.

\item For any finitely generated differential $K$-algebra $S$ with no zero divisors,
there is an algebra homomorphism $\psi: S \ra K$.

\item  For any finitely generated differential $K$-algebra $S$
there is an algebra homomorphism $\psi: S \ra K$.
\een
\ethm

\bpr
By the previous theorem, (1) implies (2).  That (2) implies (1) can be established
exactly as in the proof of Theorem~\ref{thm:3}.  Likewise for the equivalence of 
(2), (3), and (4).
\epr

\bco
A differential field $K$ in characteristic $p$ is differentially closed if and only if 
every differential variety over $K$ has a $K$-rational point.
\eco

\bpr
This is a restatement of the equivalence of (1) and (3) from the previous theorem.
\epr

We also obtain a characteristic $p$ differential Nullstellensatz, which is proved 
in exactly the same way as Corollary~\ref{diff.null}.

\bco
[Differential Nullstellensatz]
Let $K$ be a differentially closed field in characteristic $p$ with commuting derivations.
Let $\mc{I} \sub K \{x_1, \ldots , x_n  \}$ be a non-trivial differential ideal.
There is an $n$-tuple $\abar \in K$ such that for all polynomials $f(\xbar) \in \mc{I}$,
$f(\abar) = 0$.
\eco

We conclude this section with two questions.

\bqu
Let $f: X \ra Y$ be a morphism of differential finite type between affine differential  schemes
over a differential field of characteristic $p$.
Suppose that $Y$ is noetherian.  Is $f(X)$ a constructible subset of $Y$?
\equ

\bqu
\label{quest:p.sep}
Fix a prime $p \ne 0$.  Let $S$ be a differential ring of characteristic $p$, with no zero divisors,
and let $R$ be a differential subring over which $S$ is finitely generated and separable.
Let $b \in S$ be nonzero.  Is there a nonzero $a  \in R$ such that any homomorphism,
$\phi: R \ra K$ with $\phi(a) \ne 0$ lifts to a homomorphism $\psi: S \ra K$ 
with $\psi(b) \ne 0$?
\equ

\section{Difference fields}
\label{sec:ACFA}

A difference ring is a ring $R$ equipped with an injective endomorphism $\sigma$.
The classic reference on difference algebra is Cohn's book~\cite{Coh}.
Levin's book~\cite{Lev} provides an updated treatment of the subject.
Model theorists have introduced the notion of a difference closed fields,
which has become an active area of research~\cite{Mac, CH, H-ETFA}.
In this brief section, we observe that the analog of Theorem~\ref{thm:main}
fails for difference closed fields.  In a separate paper~\cite{Ros}, we prove
the difference analog of Theorem~\ref{thm:main.p}, which implies a difference
Nullstellensatz.

\bdf
Let $(K,\sigma)$ be a difference field, and let $x_1, \ldots, x_n$ be a 
set of indeterminates.  The {\em difference polynomial ring}
$K[x_1, \ldots , x_n]_\sigma$ is the polynomial ring 
$$
K[x_1, \ldots , x_n, x_1^\sigma, \ldots , x_n^\sigma, x_1^{\sigma^2}, \ldots , x_n^{\sigma^2}, \ldots ]
$$
in infinitely many variables, with $\sigma(x_n^{\sigma^m}) = x_n^{\sigma^{m+1}}$.
A difference polynomial $f \in K[x_1, \ldots , x_n]_\sigma$
determines a function $f:K^n \ra K$ in an obvious way.
\edf

A homomorphism of difference rings is a ring homomorphism that commutes with the 
endomorphisms.


\bdf
A difference field $(K,\sigma)$ is {\em difference closed} if any finite set of
difference polynomial equations and inequalities that has a solution in some difference extension
field of $K$ already has a solution in $K$.
\edf

The class of difference closed fields is first-order axiomatizable, but the theory does not have 
quantifier elimination~\cite{Mac}.  Model theorists call such fields `models of ACFA'.
The only fact needed below is that any difference field embeds in a difference closed
field, which is not, though, in any sense unique.

\bpro
\label{pr:acfa}
There are finitely generated difference rings $R \sub S$ and a difference closed field
$K$ such that there is an injective homomorphism $\phi: R \ra K$
that does not lift to any homomorphism $\psi: S \ra K$.
\epro

\bpr
Let $S = \mathbb{Q}[x]$ endowed with the automorphism $\sigma_S(x) =  -x$,
and let $R = \mathbb{Q}[x^2]$.  Note that $\sigma$ acts trivially on $R$.
Let $F$ the difference field $\mathbb{Q}(y)$, endowed with the trivial automorphism
$\sigma_F$, and let $K$ be any difference closed field containing $F$.
Then the homomorphism $\phi: R \ra F$, $x\mapsto y$, does not lift to $S$.
\epr

We now show that the analog of Theorem~\ref{thm:geo} fails for difference schemes.
Difference schemes were introduced by Hrushovski in~\cite{H-ETFA}.  We 
briefly recall the necessary background.  In order to simplify the presentation,
we assume that rings have no zero divisors.  

Let $(R, \sigma)$ be a difference ring.  A {\em transformally prime ideal} is
a prime ideal $\mfp$ such that, for all $a \in R$, $a \in \mfp$ if and only if 
$\sigma(a) \in \mfp$.  The {\em difference spectrum}, $\Spec^\sigma(R)$,
is a topological space whose elements are the transformally prime ideals.  
For each ideal $I \sub R$, the set of $\mfp \in \Spec^\sigma(R)$ containing 
$I$ is a closed set, and all closed sets arise in this way.  The space $\Spec^\sigma(R)$
can be endowed with a sheaf of difference rings in a natural way, making it into an 
{\em affine difference scheme}.  One can define abstract difference schemes
and morphisms between them, but for our purposes it suffices to know that
a homomorphism $\phi: R \ra S$ of difference rings determines a morphism
$f: \Spec^\sigma(S) \ra \Spec^\sigma(R)$.  Say that $f$ is of finite type if
$S$ is finitely generated over $R$ as a difference ring.  We will also use
the fact that if $R$ is finitely generated over a difference field, then
$\Spec^\sigma(R)$ is noetherian (\cite{H-ETFA}, Remark~3.1).

The following proposition is due to Hrushovski.

\bpro
\label{pro:Hr}
There exists a pair $X,Y$ of noetherian affine difference schemes, and a morphism
$f: X \ra Y$, such that $f(X)$ is not  constructible.
\epro

\bpr
Let $R =( \QQ[x], \sigma_R)$, with $\sigma_R(x) = x$, and let $S = (\QQ[x], \sigma_S)$,
with $\sigma_S(x) = -x$.  Define $X = \Spec^\sigma(S)$, $Y = \Spec^\sigma(R)$,
and let $f: X \ra Y$ be the morphism induced by the homomorphism 
$\phi: R \ra S$, $f(x) = x^2$.  In particular, given a transformally prime ideal 
$\mfp \sub S$, $f(\mfp) = \phi^{-1}(\mfp)$.

As $\sigma_R$ is trivial, $Y = \Spec^\sigma(R) = \Spec(R)$ and for each $n \in \QQ$,
the ideal $(x-n)$ is in $Y$.  We claim that for nonzero $n \in \QQ$, $f^{-1}(x-n) = \emptyset$
if and only if $n$ is a square.  In one direction, if $n$ is not a square, then 
$(x^2 - n) \in f^{-1}((x-n))$.  In the other direction, let $n$ be a nonzero square and suppose
for contradiction that $\mfp \in f^{-1}((x-n))$.  Then $x^2 - n \in \mfp$, but as $(x^2 - n)$
is not a prime ideal, there must be some element $a + bx \in \mfp$, $a,b \in \QQ$.  
As $\mfp$ is transformally prime $\sigma(a + bx) = a - bx$ is in $\mfp$, as is
$(a + bx) + (a -bx) = 2a$, which is invertible, unless $a = 0$.  Thus $bx \in \mfp$,
but this implies that $x \in \mfp$ and $n \in \mfp$, which is impossible.

To prove that $f(X)$ is not constructible, it now suffices to prove that
$f(X)$ is not contained in any proper closed subset of $Y$, and that 
$f(X)$ does not contain any non-empty open subset of $Y$.  First, suppose
for contradiction that there is some nonzero ideal $I \sub \QQ[x]$, with
$I \sub (x- n)$, for each nonsquare $n \in \QQ$.  Then for any polynomial $g \in I$,
each $(x -n)$ divides $g$, which is impossible.  So $f(X)$ is not contained in a 
proper closed subset of $Y$.  Arguing in the other direction, suppose again 
for contradiction that there is a nonempty $U \sub Y$ with $U \sub f(X)$.
Then for each nonzero square $n$, the ideal $(x-n)$ is in $Y \setminus U$.
But, as in the previous argument, this cannot happen.
\epr

The example in the above argument can be adapted to establish the failure of the difference
analog of Kolchin's lifting theorem, Corollary~\ref{cor:kol}.

\bpro
There is a difference ring $S$ with no zero divisors, and a subring $R \sub S$ over which
$S$ is finitely generated, with the following property.  For any nonzero $a \in R$, there is 
a transformally prime ideal $\mc{P} \sub R$, with $a \not\in \mc{P}$ , such that there does 
not exist a transformally prime ideal $\mc{Q} \sub S$, with $\mc{Q} \cap R = \mc{P}$.
\epro

\bpr
Let $S = \QQ[x]$, with $\sigma(x) = -x$, and let $R = \QQ[x^2]$, so $\sigma$ is trivial on $R$.  
(The embedding $R \ra S$ here is equivalent to the homomorphism $\phi: R \ra S$, from
the previous proof.)
For any nonzero polynomial $a \in R$, there is a nonzero square $n \in \QQ$ such that 
$x^2 - n$ does not divide $a$.  Letting $\mc{P} = (x^2-n) \sub R$, which is a transformally prime ideal,
by the proof of the preceding proposition there is no transformally prime ideal $\mc{Q} \sub S$,
with $\mc{Q} \cap R = \mc{P}$.
\epr

\section{Further results}

\subsection{Chevalley's theorem from quantifier elimination}

In this section, we give a proof of Chevalley's original theorem using
Tarski's quantifier elimination for algebraically closed fields, which can
be formulated as follows.  (Compare also Lemma~\ref{lm:qe}.)

\bthm
[Tarski]
Let $K$ be an algebraically closed field, and let  $p_1(\xbar, y), \ldots , p_j(\xbar,y),
q(\xbar,y)$ in $\ZZ[\xbar, y]$ be a finite set of polynomials. 
Suppose that there is an tuple $\cbar \in K$, and $t \in K$, such that for all 
$i \leq j$, $p_i(\cbar, t) = 0$, and $q(\cbar, t) \neq 0$.  

Then there are polynomials 
$f_1(\xbar), \ldots ,f_k(\xbar), g(\xbar) \in \mathbb{Z}[\xbar]$, such that:
\ben
\item for all $i \leq k$, $f_i(\cbar) = 0$, and also $g(\cbar) \neq 0$;
\item for any algebraically closed field $L$, and any $n$-tuple $\ebar \in L$, if 
for all $i \leq k$, $f_i(\ebar) = 0$, and also $g(\ebar) \neq 0$, then there is a $v \in K$,
such that for all $i \leq j$, $p_i(\ebar,v) = 0$, and $q(\ebar,v) \neq 0$.
\een
\ethm

\bthm[Chevalley]
Let $R \sub S$ be integral domains, such that $S$ is finitely generated over $R$
and $R$ is noetherian.
For any nonzero $b \in S$, there is a nonzero $a \in R$ with the following property.
Any homomorphism $\phi: R \ra K$, $K$ an algebraically closed field, with 
$\phi(a) \neq 0$, can be lifted to a homomorphism $\psi: S \ra K$, with $\psi(b) \neq 0$.
\ethm

\bpr
As in the proof of Theorem~\ref{thm:main}, we can reduce to the case where $S$ is
generated by a single element over $R$, $S = R(t)$, and $b = g(t)$, for some polynomial
$g(x) \in R[x]$.  Any homomorphism $\phi: R \ra K$
extends in a natural way to a homomorphism of polynomial rings $\Phi:R[x] \ra K[x]$.
Given $f(x) \in R[x]$, we write $f^\phi(x)$ for $\Phi(f(x))$.
If $t$ is transcendental over $R$, choose $a \in R$ to be the leading coefficient of $g$.
For any homomorphism $\phi: R \ra K$ with $\phi(a) \neq 0$,
$g^\phi(x) \neq 0$.  Choose $v \in K$ with $g^\phi(v) \neq 0$, and
let $\psi: S \ra K$ be the unique homomorphism lifting $\phi$ with $\psi(t) = v$.

So we may suppose that $t$ is algebraic over $R$.  Embed $S$ in an algebraically closed
field $F$.  Let $\mfp \sub R[x]$ be the prime ideal consisting of those polynomials 
$f$ such that $f(t) = 0$.  By Hilbert's basis theorem, $R[x]$ is noetherian, so let
$f_1, \ldots , f_j$ be a finite set of generators of $\mfp$.
Let $\cbar \in R$ be an enumeration of the coefficients in all the $f_i$.  Thus we can write
$f_i(x) = \hat{f}_i(\cbar,x)$, with $\hat{f}_i(\ybar, x) \in \ZZ[\ybar,x]$.  By Tarski, there is a set of 
polynomials $p_i(\ybar)$, $i \leq k$, and $q(\ybar)$, all in $\ZZ[\ybar]$,
such that 
\ben
\item  $p_i(\cbar) = 0$, for all $i$, and $q(\cbar) \neq 0$.

\item  In any algebraically closed field $L$, for any $\ebar$, if 
$p_j(\ebar) = 0$, for all $j$, and $q(\ebar) \neq 0$, then there is a $m$ such that
$\hat{f}_i(\ebar,m) =0$, for all $i$, and $\hat{g}(\ebar,m) \ne 0$.
\een
Define $a \in R$ to be $q(\cbar)$.

Let $\phi: R \ra K$ be a homomorphism with $\phi(a) \neq 0$.  Since we have
$q(\phi(\cbar)) = \phi(a)$, and $p_i(\phi(\cbar)) = 0$, for all $i$, there is a
$v \in K$ such that $\hat{f}_i(\phi(\cbar),v) = f^\phi_i(v) = 0$, for all $i$, and $g(\phi(\cbar),v) \neq 0$.
In order to establish that there is a homomorphism $\psi: S \ra K$ lifting 
$\phi$, with $\psi(t) = v$, it now suffices to show that for all $h(x) \in R[x]$,
if $h(t) = 0$, then $h^\phi(v) = 0$.  Any such $h(x)$ is equal to 
$\sum_if_i(x)w_i(x)$, with $w_i(x) \in R[x]$.  Then
$$
h^\phi(v) = \sum_if^\phi_i(v)w^\phi_i(v) =  0
$$
as desired.
\epr

\subsection{Chevalley and quantifier elimination}
\label{sec:Chev}
Given the connection between Chevalley's homomorphism extension theorem 
and quantifier elimination, one may ask for an abstract model-theoretic version
of the theorem.  
We pose a natural question in this direction, but then give an easy example
to show the answer is no.
We then formulate a restricted version of the question.
For model-theoretic terminology, see, e.g.,~\cite{Hod} or \cite{Mar}.

\bdf
Let $\mc{K}$ be a class of structures.  Say that $\mc{K}$ has the
{\em strong Chevalley property} if the following conditions hold.
Let $\mathcal{B}$ be a structure that can be embedded
in some $\mathcal{N} \in \mathcal{K}$, and let $\mathcal{A} \sub \mathcal{B}$
be a substructure such that $\mathcal{B}$ is finitely generated over $\mathcal{A}$.
Let $\bbar \in B$ be a tuple and let  $\theta(\xbar)$ be a quantifier free formula 
such that $\mathcal{B} \models \theta(\bbar)$.
Then there is a tuple $\abar \in A$ and a
quantifier free formula $\eta(\ybar)$, with $\mathcal{A} \models \eta(\abar)$,
that have the following property.  
For any $\mathcal{M} \in \mathcal{K}$, any homomorphism 
$f: \mathcal{A} \ra \mathcal{M}$ with $\mathcal{M} \models \eta(f(\abar))$
can be extended to a homomorphism $g: \mathcal{B} \ra \mathcal{M}$
 such that $\mathcal{M} \models \theta(g(\bbar))$.
\edf

\begin{question}
\label{q:chev}
Does every class of models $\mathcal{K}$ of a complete first-order theory with 
quantifier elimination have the Chevalley property?
\end{question}

\bpro
The answer to the preceding question is no.
\epro

\bpr
Let $\mc{K}$ be the class of infinite complete graphs, that is, structures $\mc{M}$ endowed
with a single binary relation $E$ such that for all $a,b \in M$, $\mc{M} \models Eab$ if any only if 
$a \neq b$. Let $\mc{M}$ be any structure in $\mc{K}$, let $\mc{A} = \mc{M}$, and let
$\mc{B} \in \mc{K}$ extend $\mc{A}$ by a single element.
Then the identity homomorphism $f: \mc{A} \ra \mc{M}$ does not lift to $\mc{B}$.

Alternately, let $\mc{K}$ be the class of dense linear orders without endpoints, let
$\mc{A} = \mc{M} = \mathbb{Q}$ and let $\mc{B} = \mathbb{Q} \cup \{\sqrt{2}\}$, all
with the canonical order.
\epr

(We note that it is also easy to construct a counterexample where the 
language contains only function symbols.)

One might also ask whether a class $\mc{K}$ not having quantifier elimination
implies that it does not have the strong Chevalley property.
An affirmative answer would show that Propositions~\ref{pr:dcfp} and~\ref{pr:acfa}
follow immediately from the failure of quantifier elimination for
$DCF_p$ and $ACFA$.
Nevertheless, we observe that this is not the case.  Thus
quantifier elimination is neither necessary nor sufficient for a class to have 
the strong Chevalley property.

\bpro
There is a complete first-order theory $T$ without quantifier
elimination, such that the class $\mc{K}$ of models of $T$ has the
strong Chevalley property.
\epro

\bpr
Let $L$ contain a single binary relation symbol $E$, and let $T$ be the 
complete theory that says that $E$ is an equivalence relation, every
equivalence class has size 2 or 3, and there are infinitely many 
equivalence classes both of size 2 and of size 3.

Suppose now that $\mc{B}$ is a substructure of some model of $T$,
and let $\mc{A} \sub \mc{B}$ be such that $\mc{B} \setminus \mc{A}$
is finite.  Given $\bbar \in B$ and a quantifier free $\theta(\xbar)$
such that $\mc{B} \models \theta(\bbar)$, let $\abar \in A$
be the set of all $a \in A$ that are $E$-equivalent to some $b \in \bbar$,
and let $\eta(\ybar)$ be the complete atomic diagram of $\abar$.
It is easy to verify that $\abar$ and $\eta(\ybar)$ have the
desired properties.
\epr

On the other hand, there are well-known model theoretic facts related to 
the above question, where homomorphisms are replaced
by embeddings.  Here is an example.

\bpro
Let $\mc{K}$ be the class of models of a complete first-order theory with
quantifier elimination.  Let $\lambda$ be an infinite cardinal.  
Let $\mc{B}$ be a structure of size $\leq \lambda$ that can be embedded in 
some $\mc{N} \in \mc{K}$, and let $\mc{A}$ be a substructure of
size $ < \lambda$.  For any $\lambda$-saturated model $\mc{M} \in \mc{K}$,
any embedding $f: \mc{A} \ra \mc{M}$ extends to an embedding
$g: \mc{B} \ra \mc{M}$.
\epro

We now formulate a refined version of Question~\ref{q:chev}.

\bdf
Let $\mc{K}$ be a class of structures.  Say that $\mc{K}$ has the
{\em Chevalley property} if the following conditions hold.
Let $\mathcal{B}$ be a finitely generated structure that can be embedded
in some $\mathcal{N} \in \mathcal{K}$, and let $\mathcal{A} \sub \mathcal{B}$
be a finitely generated substructure.
Let $\bbar \in B$ be a tuple and let  $\theta(\xbar)$ be a quantifier free formula 
such that $\mathcal{B} \models \theta(\bbar)$.
Then there is a tuple $\abar \in A$ and a
quantifier free formula $\eta(\ybar)$, with $\mathcal{A} \models \eta(\abar)$,
that have the following property.  
For any $\mathcal{M} \in \mathcal{K}$, any homomorphism 
$f: \mathcal{A} \ra \mathcal{M}$ with $\mathcal{M} \models \eta(f(\abar))$
can be extended to a homomorphism $g: \mathcal{B} \ra \mathcal{M}$
 such that $\mathcal{M} \models \theta(g(\bbar))$.
\edf

\begin{question}
Does every class of models $\mathcal{K}$ of a complete first-order theory with 
quantifier elimination have the Chevalley property?
\end{question}

\providecommand{\bysame}{\leavevmode\hbox to3em{\hrulefill}\thinspace}
\providecommand{\MR}{\relax\ifhmode\unskip\space\fi MR }
\providecommand{\MRhref}[2]{%
  \href{http://www.ams.org/mathscinet-getitem?mr=#1}{#2}
}
\providecommand{\href}[2]{#2}


\end{document}